\documentclass[12pt,leqno]{amsart}	
\usepackage{amsthm,amsmath,amsfonts,amscd,verbatim,latexsym,exscale,xypic}
\usepackage[mathscr]{eucal}
\input{defns.cr}

\begin{document}
\title{On Equations Defining Coincident Root Loci}
\author{Jaydeep V. Chipalkatti}
\maketitle

\bigskip 
\parbox{12cm}{
{\bf Abstract.}
We revisit an old problem in classical invariant theory, viz.~that of 
giving algebraic conditions for a binary form to have linear factors 
with assigned multiplicities. We construct a complex of 
$SL_2$--representations such that the desired algebraic 
conditions are expressible as a specific cohomology group of this complex.
}

\medskip 

\section{Introduction}
Consider a binary form $F(x,y)=\sum\limits_{j=0}^d  a_j\, x^jy^{d-j}$
over the complex numbers. It splits as a product of linear forms, 
and it is classical that $F$ has a repeated factor if and only if its 
discriminant vanishes. Similarly, we can ask for algebraic conditions on 
the coefficients $a_j$, so that $F$ has say, a triple factor or two double 
factors. More generally, we may fix a partition $\lambda$ of $d$, and 
ask for algebraic conditions so that the factors of $F$ have multiplicities 
as dictated by the parts in $\lambda$. The object of this note is to devise a 
method for answering such questions.

This problem is addressed for the first time (to my knowledge) by 
Arthur Cayley \cite{Cayley1}. Thereafter J.~Weyman has obtained substantial 
results for special kinds of partitions 
in \cite{Weyman1}, \cite{Weyman3} and \cite{Weyman2} (see 
section \ref{comments.section}).

The polynomial $F(x,y)$ as above will be identified with the point
$[a_0,\dots,a_d]$ of $\P^d$. Let then 
$\lambda=(\lambda_1,\dots,\lambda_n)$ be a partition of $d$, and 
consider the set \[ X_\lambda = \{ [a_0,\dots,a_d]: 
 F \; \text{factors as \,$\prod\limits_{i=1}^n L_i^{\lambda_i}$\,
            for some linear forms $L_i$} \}, \]
which is a projective subvariety of $\P^d$. This is the coincident root locus 
in the title and the principal object of study in this paper. 

The group $\SL_2({\mathbb C})$ acts on the 
imbedding $X_\lambda \subseteq \P^d$ in the following manner. 
An element $ \left( \begin{array}{cc} 
\alpha & \beta \\ \gamma & \delta \end{array} \right) \in \SL_2({\mathbb C})$ 
sends the form $F$ to $F'=\sum\limits_{j=0}^d  a_j'\, x^jy^{d-j}$,
where the $a_j'$ are determined by 
\[ \sum\limits_{j=0}^d  a_j\, 
(\alpha\, x + \beta\, y)^j(\gamma\, x +\delta\, y)^{d-j}= 
\sum\limits_{j=0}^d {a_j'}\, x^jy^{d-j}. \] 
Its defining ideal $\,I_{X_\lambda} < {\mathbb C}[a_0,\dots,a_d]$\, 
is then an $\SL_2$--subrepresentation, and our problem is essentially 
one of calculating this ideal. Our main result (Corollary 
\ref{corollary.gradedpiece})
is a (somewhat indirect) solution to this problem.  
Specifically, for each positive integer $m$, we construct a complex of 
$SL_2({\mathbb C})$--representations (see formula (\ref{qalphap.formula}))
whose zeroth cohomology group is the graded piece $(I_X)_m$. 
In section \ref{computation.section}, we use this complex to 
calculate the ideal for a few specific examples and express the answer 
in the language of classical invariant theory (especially see 
Theorems \ref{criterion.32} and \ref{criterion.33}). In 
Theorem \ref{sing.theorem} of section \ref{singularities.section},
we identify the singular locus of $X_\lambda$; 
a result which should be of interest in itself. 

The base field will be ${\mathbb C}$. All terminology from algebraic 
geometry follows Hartshorne \cite{Ha}. The formulae are numbered 
in a single sequence throughout the text. 

\bigskip
\noindent {\sc Acknowledgements:\,} 
The program Macaulay--2 has been of immense help in computations, and it is 
a pleasure to thank its authors Dan Grayson and Mike Stillman. 
The financial assistance of Anthony V. Geramita, Leslie Roberts and 
Queen's University during the progress of this work 
is gratefully acknowledged.

\section{Preliminaries} 
In the next three subsections, we recall a few matters from the invariant 
theory of binary forms. An excellent reference for the classical 
theory is the text by Grace and Young \cite{GrYo}. See 
\cite{FH} and \cite{Springer1} for the modern theory and 
\cite{Sturmfels} for a discussion of algorithms for the computation 
of invariants and covariants. 

\subsection{Representations of $SL_2({\mathbb C})$.}
In the sequel, $V$ denotes a two dimensional ${\mathbb C}$--vector space.
We will use the following isomorphisms of 
$SL(V)$--modules\,:
\begin{equation} \begin{array}{ll}
\wedge^m(\Sym^n \,V) & = \Sym^m(\Sym^{n+1-m}\, V),   \\
\Sym^m \,V \otimes \Sym^n \,V & =   \opluslim\limits_{r=0}^{[\frac{m+n}{2}]}
\Sym^{m+n-2r} \,V \quad \text{(Clebsch--Gordan formula);} 
\label{slv.identities} \end{array} \end{equation}
and the Cayley--Sylvester formula 
\begin{equation} 
\Sym^m(\Sym^n \,V)=  \opluslim\limits_{r=0}^{[\frac{mn}{2}]}
(\Sym^{mn-2r} \,V )^{\,\oplus\, p(r,m,n)-p(r-1,m,n)}, 
\label{plethysm} \end{equation}
where $p(r,m,n)$ is the number of partitions of 
$r$ into at most $m$ parts with no part exceeding $n$. The $GL(V)$--modules 
$\Sym^d \, (V^*)$ and $(\Sym^d \, V)^*$ are isomorphic, and we write them
indifferently as $\Sym^d \, V^*$. 

Let $\{x,y\}$ be a basis of $V$. With the identification 
$a_j=(x^jy^{d-j})^* \in \Sym^d\, V^*$, 
the form $F=\sum a_j \,x^jy^{d-j}$ corresponds to the trace element 
$\sum\limits_j (x^jy^{d-j})^* \otimes x^jy^{d-j}$ in 
$\Sym^d \,V^* \otimes \Sym^d \,V$.
The ring $\mathbb C[a_0,\dots,a_d]$ is identified with the symmetric 
algebra $\Sym^\bullet\, (\Sym^d\,V^*)$. 
\subsection{Invariants and Covariants} \label{cov.subsection}
Let 
\[ \Sym ^q \,V^* \hookrightarrow  \Sym^p(\Sym^d \,V^*)\] 
be an $SL(V)$--subrepresentation; or what is the same, a nonzero 
$SL(V)$--morphism 
${\mathbb C} \lra \Sym^p(\Sym^d \,V^*) \otimes \Sym^q \,V$. Its 
image consists of all scalar multiples 
of a polynomial $\Phi(a_0,\dots,a_d;x,y)$, having total degree $p$ in the 
$a_j$ and $q$  in $x,y$. This polynomial is called a 
covariant of degree $p$ and order $q$ (of the binary form $F$). We will also 
say that $\Phi$ is of type $(p,q)$. 
E.g., for $d \ge 2$, the inclusion 
$\Sym^{2d-4} \,V^* \hookrightarrow \Sym^2(\Sym^d \,V^*)$
defines a covariant of type $(2,2d-4)$, called the Hessian of $F$. 

An invariant is a covariant of order zero. E.g., the discriminant 
of $F$ (defined to be the Sylvester resultant of 
${\partial F}/{\partial x}$ and ${\partial F}/{\partial y}$) is 
an invariant of degree $2d-2$.

\subsection{Transvectants.}
Let $\Phi_1, \Phi_2$ be covariants of types $(p_1,q_1),(p_2,q_2)$ 
respectively, and $1 \le r \le \min \{q_1,q_2 \}$ an integer. 
Associated to these is a covariant of type $(p_1+p_2,q_1+q_2-2r)$ 
called the $r$-th transvectant of $\Phi_1, \Phi_2$, written 
$(\Phi_1,\Phi_2)^r$. (By convention
$(\Phi_1,\Phi_2)$ stands for $(\Phi_1,\Phi_2)^1$.) 
For instance, the Hessian equals $(F,F)^2$. 
\begin{Remark} \rm \label{transvectant.remark}
The transvectant can be defined in several ways, and the various 
definitions are in agreement only upto a scalar. In 
section \ref{computation.section}, it will be necessary for us 
to fix this scalar, and it is important that some 
definite convention be followed consistently.  

We will use the definition given in terms of the symbolic method.
Let\footnote{The introduction of binomial coefficients is a familiar 
device in order to make things work smoother.}
$F=\sum\limits_{j=0}^d  \binom{d}{j} a_j \,x^jy^{d-j}$, and 
\[ \Phi_i= \sum\limits_{k=0}^{q_i} \varphi_{ik} \, x^k y^{q_i-k}, 
\quad i=1,2. \] 
Thus each $\varphi_{ik}$ is a homogeneous form of degree 
$p_i$ in the $a_j$.

Now formally write 
$ \Phi_1=(\alpha_1 \,x+ \alpha_2 \,y)^{q_1},\; 
   \Phi_2=(\beta_1 \,x + \beta_2 \,y)^{q_2} $ and 
define 
\[ (\Phi_1,\Phi_2)^r = (\alpha_1 \,\beta_2 -\alpha_2 \,\beta_1)^r 
(\alpha_1 \,x+ \alpha_2 \,y)^{q_1-r}(\beta_1 \,x +\beta_2 \,y)^{q_2-r}.
\]
That is to say, expand the right hand side and substitute 
$\varphi_{1k}$ and $\varphi_{2k}$ respectively, for
$\binom{q_1}{k} \alpha_1^k\,\alpha_2^{q_1-k}$ and 
$\binom{q_2}{k} \beta_1^k\,\beta_2^{q_2-k}$.
The symbolic method is more sound than it appears prima facie,
see \cite{KungRota} for a thorough modern treatment. 
\end{Remark}
\subsection{Definition of the Coincident Root Locus} 
\label{crlocus.subsection}
\setcounter{footnote}{1}

In the sequel $\lambda=(1^{e_1}2^{\,e_2} \dots d^{\,e_d})$ denotes a partition 
of $d$, having $e_r$ parts of size $r$ for $1 \le r \le d$. 
The number of parts is $\sum e_r =n$. 

With $\{x,y\}$ a basis of $V$, we identify a 
point of $\P\, \Sym^{e_r} V$ with a degree $e_r$ polynomial 
$G_r(x,y)$ determined upto scalars.

Now we have Veronese maps 
\begin{equation}
 v_r: \P\, \Sym^{e_r} \,V \lra \P\, \Sym^{r \,e_r} \,V, 
\quad G_r \lra G_r{}^r \end{equation} and a multiplication map 
\begin{equation}
\mu: \prod\limits_{r=1}^d  \P\, \Sym^{r\,e_r}\,V \lra \P\, \Sym^d \,V, 
\quad (H_1,\dots,H_d) \lra \prod\limits_{r=1}^d H_r. 
\label{multiplication} \end{equation}
Finally consider the composite $f_\lambda$,
\begin{equation} \label{defn.nulambda} 
\begin{aligned}
(Y_\lambda =) & \prod_r \P\, \Sym^{e_r} \,V \stackrel{\Pi v_r} \lra \prod_r 
\P\, \Sym^{r \,e_r} \,V \stackrel{\mu}{\lra} \P\, \Sym^d \,V, \\
  & (G_1,\dots,G_d) \stackrel{f_\lambda}{\lra} \prod\limits_{r=1}^d G_r{}^r.
\end{aligned} \end{equation}
\begin{Definition} \rm 
The Coincident Root locus $X_\lambda$ is defined to be the scheme--theoretic 
image of $f_\lambda$. \end{Definition} 
The locus is closed in $\P\, \Sym^d\, V$ since $f_\lambda$ is a projective 
morphism, and it is reduced since $Y_\lambda$ is. Its dimension is 
the number of parts in $\lambda$, which is $n$. 

The morphism $Y_\lambda \stackrel{f_\lambda}{\lra} X_\lambda$ is finite
(being quasi-finite and projective) and birational (since for a general $F$
in the image, $G_r$ can be recovered as the product of linear forms
dividing $F$ exactly $r$ times). The ideal sheaf $\I_{X_\lambda}$ is 
the kernel
\begin{equation} 0 \lra \I_{X_\lambda} \lra \O_{\P^d} \lra 
\O_{X_\lambda} \lra 0, \label{defn.ix} \end{equation}
and for later use we define $\D_\lambda$ to be the cokernel
\begin{equation}
 0 \lra \O_{X_\lambda} \lra {f_\lambda}_* \O_{Y_\lambda} \lra \D_\lambda 
\lra 0. \label{defn.dlambda} \end{equation}
We will drop the suffix $\lambda$ when it is safe to do so. 
\begin{Remarks} \rm 
\begin{itemize}
\item{
We have $X_{\lambda_1} \subseteq X_{\lambda_2}$ exactly when $\lambda_2$ is 
a refinement of $\lambda_1$.}
\item{
The locus $X_{(d)}$ is the rational normal curve of degree $d$ and 
$X_{(d-1,1)}$ its tangential developable. In general $X_{(d-s,1^s)}$ is 
the tangential developable (i.e., the closure of the union of tangent spaces 
at smooth points) of $X_{(d-s+1,1^{s-1})}$. Of course,
$X_{(2,1^{d-2})}$ is the hypersurface of degree $2(d-1)$ defined by 
the discriminant.}
\item{ The degree of the CR locus in $\P^d$ is given by 
\[ 
\deg X_\lambda = 
\frac{n!}{\prod\limits_r (e_r !)} \prod\limits_r r^{e_r}, \] 
a formula going back to Hilbert. To see this, identify 
$\P\,\Sym^d \,V \cong \Sym^d\,(\P^1)$ with 
the set of effective divisors of degree $d$ on $\P^1$. Let $\Sigma$ denote 
a general $(d-n)$-dimensional linear subspace of $\P^d$, it corresponds to 
a linear series $g^d_{d-n}$ on $\P^1$. The points of $\Sigma \cap X$ correspond 
to those divisors in the series which may be written as 
$\sum\limits_r r(P_{r,1}+P_{r,2}+ \dots+P_{r,e_r})$, for some points 
$P_{r,j} \in \P^1$. According to De Jonqui{\`e}res' 
formula (see \cite[p. 359]{ACGH}), the number of such divisors is 
the coefficient of \,$t_1^{\,e_1}t_2^{\,e_2}\dots t_d^{\,e_d}$ in 
$(1+t_1+2\,t_2+\dots+d\,t_d)^n$, hence the assertion.}
\end{itemize} \end{Remarks}

\section{The ideal of the CR locus}  \label{crideal.section}
\subsection{The Machine Computation of $I_X$.} 
There is a rather straightforward algorithm for the computation of 
the ideal $I_X$, which can be implemented on any computer algebra 
system with modest capabilities. Since $X$ is given as the image of 
a projective morphism, the problem is one of elimination 
theory. We illustrate with the case $\lambda= (3\,2\,2)$. 
\label{machine.subsection}
Consider the morphism 
\[ f_{322}: \P\, \Sym^2 \,V \times \P\, V \lra 
   \P\, \Sym^7 \,V, \quad (G_2,G_3) \lra G_2^2\,G_3^3\,(=F). \] 
Using affine c{\"o}ordinates, write 
\[ 
   G_2 = t^2+u_1\,t+u_2, \; G_3=t+v_1, \; 
   F = t^7+ a_1\,t^6+ \dots + a_7. \] 
By forcing the equality 
\[ t^7+ a_1\,t^6+ a_2\,t^5+ \dots + a_7
=(t^2+u_1\,t + u_2)^2(t + v_1)^3,  \] we get polynomial expressions 
$ a_j=q_j(u_1,u_2,v_1)$ for $1 \le j \le 7$.
This defines a homomorphism 
\[ q_{322}: {\mathbb C}[a_1,\dots,a_7] \lra
{\mathbb C}[u_1,u_2,v_1], \] which is the ring theoretic counterpart 
of $f_{322}$. Calculate a Gr{\"o}bner basis for 
the kernel of $q_{322}$ (with respect to a degree lex or degree 
reverse--lex term order) 
and homogenize with respect to a new variable $a_0$; 
this gives the ideal $I_X$. (See e.g. \cite[Exer. 1.6.19c]{AdLou}.) 

The actual calculation shows that it is generated by a 
$364$-dimensional subspace of ${\mathbb C}[a_0,\dots,a_7]_6$. 
This is not very informative by itself; and to gain some insight 
into the structure of $I_X$, we aim to describe it 
invariant--theoretically. In the present example,
this means identifying the subrepresentation $(I_X)_6 \subseteq 
\Sym^6(\Sym^7\, V^*)$. To this end, we construct our main 
technical artifice (the complex $\G^\bullet$) in the next subsection. 

The idea, in a nutshell, is to approach $X$ indirectly by working with its 
normalisation $Y$. We will construct a locally free resolution of the 
structure sheaf $\O_Y$, and then generate a spectral sequence by taking 
its pushforward.

\subsection{The Eagon--Northcott Complex} \label{encomplex.subsection}

Set $T = Y  \times \P\,\Sym^d \,V$ with projections 
$\pi_r: T \ra \P\, \Sym^{e_r} \,V$ for $1 \le r \le d$, 
and $\pi: T \ra \P\, \Sym^d \,V$. Factor the map $f$ as 
\[ Y \stackrel{1 \times f} \lra T \stackrel{\pi}{\lra} 
\P\,\Sym^d \,V,  \]
and let $\Gamma \subseteq T$ denote the image of $ 1 \times f$.
\[ \diagram
T \rto^(0.4){\pi_r}^{} \dto^{\pi} & 
\P\, \Sym^{e_r} V & Y \rto^{\text{iso}} \drto_f & 
\; \Gamma \dto^{\pi} \rto|<<\tip^{\text{incl}} & T \dto^{\pi} \\ 
\P\, \Sym^d \,V & & & X \rto|<<\tip^{\text{incl}} & { \P^d} \enddiagram \]

We claim that the structure sheaf $\O_\Gamma$ is resolved by an 
Eagon--Northcott complex of vector bundles on $T$.

Quite generally, let $\aa \stackrel{\varphi}{\lra} \bb$ be a morphism of 
vector bundles on a regular scheme $T$. Assume that the bundles have ranks 
$a,b$ respectively, with $a \ge b$. Then we have a complex 
\[ 0 \ra \N^{b-a-1} \ra \dots \ra \N^p \ra \N^{p+1} \ra \dots \ra \N^0 (=\O_T)
 \ra 0, \]
where 
\begin{equation}
 \N^p = (\wedge^b \bb)^{-1} \otimes \Sym^{-(p+1)}\, \bb^* \otimes 
\wedge^{b-p-1} \aa \quad \text{for}\;\; b-a-1 \le p \le -1.
\label{np} \end{equation}
The differential of the complex is defined via contraction with the section 
$\varphi \in H^0(T,\aa^* \otimes \bb)$. (See \cite[Appendix A2.6]{Ei} for details 
of the construction.)

Let $T_{\varphi}$ be the closed subscheme of $T$ defined by the ideal sheaf 
$\;\text{Fitt}_0(\,\coker\,\varphi)$. Set--theoretically, it is the 
{\sl degeneracy locus} $\{t \in T: \rank\; \varphi_t \le b-1 \}$.
If the subscheme has codimension $a-b+1$ and is a local complete intersection,
then $\N^{\bullet}$ is a resolution of its structure sheaf 
$\O_{T_{\varphi}}$ (loc.cit.). 

Returning to our setup, define line bundles 
\[ \ll = \otimeslim\limits_{r=1}^d \pi_r^* \O_{\P^{e_r}}(r),\;\;
\mm = \pi^* \O_{\P^d}(1)
\] on $T$. 
By the K{\"u}nneth formula, 
\[ 
H^0(T,\ll)= \otimeslim\limits_{r=1}^d \Sym^r(\Sym^{e_r} V^*), \;\;
H^0(T,\mm)= \Sym^d \, V^*. \] 
Define the {\sl comultiplication} map 
\begin{equation}\label{comultiplication}
\Sym^d \,V^* \lra \otimeslim\limits_{r=1}^d \Sym^r(\Sym^{e_r} \,V^*),
\end{equation}
as dual to the usual multiplication map on polynomials.
Adding it to the identity 
$\Sym^d \,V^* \lra H^0(T,\mm)$ gives a map 
\[ \Sym^d \,V^* \lra H^0(T,\ll \oplus \mm), \] which defines a 
morphism of bundles 
\[ \varphi: \underbrace{\Sym^d \,V^* \otimes \O_T}_{\aa} \lra 
\underbrace{\ll \oplus \mm}_{\bb}, \quad \text{with $a=d+1, b=2$.}
\] 

\begin{Theorem} 
The Eagon--Northcott complex of $\varphi$ resolves $\O_\Gamma$. 
\end{Theorem} 
\demo We will show that the locus $ \{\,\rank\; \varphi \le 1 \}$ coincides 
with $\Gamma$. Since $\Gamma$ (being smooth) is a local complete intersection 
of codimension $d$ in $T$, the theorem will follow. 

It is easier to work with the transpose map 
\[ \varphi^{\text{trans}}: \ll^{-1} \oplus \mm ^{-1} 
\lra \Sym^d\, V \otimes \O_T, \] 
which has the same degeneracy locus. Let then 
$t=(G_1,\dots,G_d; F)$ be a point of $T$. The fibre of the line 
bundle $\O_{\P^{e_r}}(-r)$ over the point $G_r \in \P\, \Sym^{e_r}V$
is naturally isomorphic to the linear span of $G_r^r$, hence the 
fibre of $\ll^{-1}$ over $t$ is the space 
${\mathbb C}\,\langle \prod\limits_r G_r^r \rangle $. So over $t$, the map 
$\varphi^{\text{trans}}$ is given by 
\[ {\mathbb C}\, \langle \prod\limits_r G_r^r \rangle
\oplus {\mathbb C}\, \langle F \rangle 
\lra {\mathbb C}\, \langle F,\prod\limits_r G_r^r \rangle. \] 
This map has rank one iff $F$ is a scalar multiple of 
$\,\prod\limits_r G_r^r$ iff $t$ lies in $\Gamma$. 
The theorem is proved. \qed 

\medskip 
After expanding (\ref{np}), 
\begin{equation} \label{npexpanded} 
\begin{aligned}
\N^p= \{ \bigoplus_{\stackrel{\alpha+\beta =-(p+1),}{\alpha,\beta \ge 0}} 
\ll^{-(\alpha+1)} \otimes \mm^{-(\beta+1)}\}
\, \otimes  \wedge^{1-p} &  (\Sym^d \,V^*) \\
& \text{for}\; -d \le p \le -1, \end{aligned} \end{equation}
\[ \N^0=\O_T \;\; \text{and} \;\;
\N^\bullet \, {\stackrel{\sim}{\ra}}_q \;\O_\Gamma. \]

Now consider the second quadrant spectral sequence 
\begin{equation} \begin{aligned}
  E_1^{p,q} & = {\mathbf R}^q\pi_*\, \N^p \quad 
\text{for $-d \le p \le 0$ and $q \ge 0$ ;} \\
  d_r^{\,p,q}: E_r^{p,q} & \lra E_r^{p+r,q-r+1}, \quad 
   E_\infty^{p,q}  \Rightarrow {\mathbf R}^{p+q}\pi_* \O_\Gamma. 
\end{aligned} \label{spectralseq1} \end{equation}
Each $E_1^{p,q}$ term can be precisely calculated using the projection 
formula followed by the K{\"u}nneth formula. 
Apart from $E_1^{0,0}=\O_{\P^d}$, all nonzero terms are in the $n$-th row.
In fact 
\[ \begin{aligned} E_1^{p,n} = 
\{  \bigoplus_{\alpha+\beta=-(p+1) \atop \alpha,\beta \ge 0} \{ 
\otimeslim_r H^{e_r}(\P^{e_r},\O_{\P}(-r \alpha -r)) \} 
\, \otimes \,  & \O_{\P^d}(-\beta-1) \} \\
\, \otimes \,  &  \wedge^{1-p}(\Sym^d \,V^*).
\end{aligned} \] 
Now $H^{e_r}(\P^{e_r},\O_{\P}(-r \alpha -r)) \neq 0$ iff 
$r \alpha + r \ge e_r+1$. Using Serre duality, we get 
\begin{equation}\begin{aligned}
E_1^{p,n} = \{  \bigoplus_{\stackrel{\alpha+\beta=-(p+1)}{\alpha \ge 
{\scriptscriptstyle M}-1,\,\beta \ge 0}}
\{ \otimeslim_r \Sym^{r\alpha +r -e_r-1}(\Sym^{e_r} \,V) \} & \otimes 
 \O_{\P^d}(-\beta-1) \}  \\ 
 & \otimes \wedge^{1-p}(\Sym^d \,V^*).  \label{rqpinp} 
\end{aligned} \end{equation}
Here ${\scriptstyle M}$ denotes $\lceil \max\limits_r \{\frac{e_r+1}{r}\} \rceil$. 
Now $\pi$ is a finite morphism on $\Gamma$, hence 
$E_\infty^{p+q}=0$ for $p+q \neq 0$ and we have an extension 
\begin{equation}
  0 \lra E_\infty^{0,0} \lra \pi_* \O_\Gamma \lra E_\infty^{-n,n} 
 \lra 0. \label{ext1} \end{equation}

\begin{Proposition} We have 
\[ E_2^{p,n} = \dots = E_{n+1}^{p,n}= \begin{cases}
\I_X \quad & \text{if $p=-(n+1)$,} \\
\D  \quad  & \text{if $p=-n$,} \\
0 \quad    & \text{otherwise.} \end{cases} \] 
\end{Proposition}
\demo All the differentials from $d_2$ upto $d_n$ are zero, hence 
$E_2^{p,n} = \dots = E_{n+1}^{p,n}$. If $p \neq -n,-(n+1)$, then 
$E_2^{p,n}= E_\infty^{p,n}=0$. Now going back to the definition 
of the filtration on the abutment of a spectral sequence 
(see e.g. \cite[Ch. III]{BT}), it follows that $E_\infty^{0,0}$ is the 
image of the morphism $\pi_* \O_T \ra \pi_* \O_\Gamma$, which is 
$\O_X$. Hence $E_2^{-n,n}= E_\infty^{-n,n}$ equals 
$\text{coker}( \O_X \ra \pi_* \O_\Gamma)= \D$. Finally 
$E_2^{-(n+1),n} = \text{ker}\, (E_1^{0,0} \ra E_\infty^{0,0}) = \I_X$.
\qed 

The proposition exhibits $\I_X$ as the middle cohomology of the complex 
\[ E_1^{-(n+2),n} \ra E_1^{-(n+1),n} \ra E_1^{-n,n}. \]
Each term of the complex is a vector bundle on $\P^d$, explicitly described 
in terms of the combinatorial datum $\lambda$. 

Our next step is to work with the complex given by the $n$-th row of this 
spectral sequence. Let $\G^\bullet$ be 
the complex $E_1^{\bullet,n}[-(n+1)]$, with differential denoted 
$d_G$. By the previous proposition, 
$\H^p(\G^\bullet)=\I_X, \D$ for $p=0,1$ respectively and zero 
elsewhere. The point of introducing the shift is merely to simplify 
the indices. $\G^p$ is nonzero in the range 
$ n+1-d \le p \le n+1- {\scriptstyle M}$. 

\subsection{Hypercohomology of $\G^\bullet$}
\label{hyperg.subsection}

Set $\G^\bullet(m) = \G^\bullet \otimes \O_{\P^d}(m)$ for $m \ge 0$.
Our interest lies in the hypercohomology of the complex 
$\G^\bullet(m)$. There are two spectral sequences 
\begin{equation} \begin{align}
E_1^{p,q} & = H^q(\P^d, \G^p(m)) , 
& d_r^{\,p,q}: E_r^{p,q} \lra E_r^{p+r,q-r+1}, \label{2specs.1} \\
{\grave E}_2^{p,q} & = H^q(\P^d, \H^p(\G^\bullet)(m)), 
& {\grave d}_r^{\,p,q}:{\grave E}_r^{p,q} \lra {\grave E}_r^{p-r+1,q+r};  
\end{align} \label{2specs.2} \end{equation} with common abutment
$\Hy^{\,p+q}(\G^\bullet(m))$. 

The term $\G^p(m)$ is a sum of line bundles of the form 
$\O_{\P^d}(m-\beta-1)$. Now $m-\beta-1 \ge m-d$ (this is clear from 
formula (\ref{rqpinp})), hence $H^{\ge 1}(\P^d,\G^p(m))=0$. 
Thus $E_1^{p,q}=0$ for $q \ge 1$, this forces $E_2=E_\infty$.
Now ${\grave E}_2^{p,q}$ is zero outside the columns $p=0,1$, so 
${\grave E}_3 = {\grave E}_\infty$. 
\begin{Theorem}  \label{hyperg.theorem} 
Assume $m \ge 0$. Then 
\[ \Hy^i = \begin{cases}
H^0(\P^d, \I_X(m)) & \text{for $i=0$,} \\
0                  & \text{for $i \neq 0,1$.}
\end{cases} \] There is an extension
\[ 0 \lra H^1(\P^d, \I_X(m)) \lra \Hy^1 \lra {\grave E}_3^{1,0} 
\lra 0, \] and 
${\grave E}_3^{1,0}$ equals the image of the map 
\[ H^0(\P^d, f_* \O_Y(m)) \lra  H^0(\P^d, \D(m)).\]
In particular, ${\grave E}_3^{1,0}= 
H^0(\P^d, \D(m))$ if $H^1(\P^d,\O_X(m))=0$.
\end{Theorem}
We will use the following technical lemma in the proof.
\begin{Lemma} For $m \ge 0$, 
the map ${\grave d}_2^{1,q}: H^q(\D(m)) \lra H^{q+2}(\I_X(m))$ is 
bijective for $q \ge 1$ and surjective for $q=0$.
\end{Lemma}
\demo 
The map is a composite $H^q(\D(m)) \ra H^{q+1}(\O_X(m)) 
\ra H^{q+2}(\I_X(m))$ of two connecting homomorphisms in the 
long exact sequences associated to (\ref{defn.ix}) and 
(\ref{defn.dlambda}). (This follows directly from the definition 
of ${\grave d}_2$ (loc.cit.).) 
Now $H^{\ge 1}(\P^d, \O_{\P^d}(m))=0$, and we claim 
that $H^{\ge 1}(\P^d, f_* \O_Y(m))=0$. Indeed, 
\[ \begin{array}{ll}
H^i(\P^d, f_* \O_Y(m)) =   H^i(Y, f^* \O_{\P^d}(m)) 
\; &\text{\;\;(by Leray spectral sequence)} \\
= \bigoplus\limits_{\sum i_r=i} \;\;
 \otimeslim\limits_{r=1}^d H^{i_r}(\P^{e_r},\O_{\P^{e_r}}(mr))
\; & \text{\;\;(by K{\"u}nneth formula)}
\end{array} \] which vanishes for $i \ge 1$. It now follows 
that the first connecting map is a bijection for $q \ge 1$ 
and a surjection for $q=0$; and the second map is always an 
isomorphism. The lemma is proved.  \qed 

\medskip 
{\sc Proof of theorem \ref{hyperg.theorem}.}

Since ${\grave E}_2^{p,q}=0$ for $p \neq 0,1$, we have $\Hy^i=0$ for 
$i < 0$. By the lemma, ${\grave E}_3^{p,q}=0$ for $p=0,q \ge 2$ and 
$p=1,q \ge 1$; hence $\Hy^i=0$ for $i \ge 2$. The rest is clear. \qed 

\begin{Corollary} \label{corollary.gradedpiece}
The space $H^0(\P^d, \I_X(m))$ of hypersurfaces of degree $m$ vanishing 
on $X$ is given by the middle cohomology of the complex 
\[ H^0(\P^d, \G^{-1}(m)) \ra 
H^0(\P^d, \G^{0}(m)) \ra H^0(\P^d, \G^{1}(m)).
\]  \qed \end{Corollary} 

At this point, the determination of this space is in principle 
a problem in linear algebra; once we have made the differential 
of the complex explicit. We will merely outline the description 
here and leave the details to the diligent reader. 

Let $d_{G,m}^{\,p}: H^0(\P^d, \G^p(m)) \lra H^0(\P^d, \G^{p+1}(m))$ 
denote the differential in question.
Define \begin{equation} \label{qalphap.formula}
 \begin{aligned}
{} & z(\alpha,r)=r\alpha +r -e_r-1,\, 
M(\alpha)= \otimeslim\limits_r \Sym^{z(\alpha,r)}(\Sym^{e_r} \,V) 
\;\; \text{and} \\ 
& Q(\alpha,p)=M(\alpha) \otimes \Sym^{m+p+\alpha-n-1}(\Sym^d \,V^*)
\otimes \wedge^{n+2-p}(\Sym^d \,V^*). \end{aligned} \end{equation}
Then 
$ H^0(\P^d, \G^p(m)) = 
\bigoplus\limits_{n-p \ge \alpha \ge {\scriptscriptstyle M}-1}
Q(\alpha,p)$, and it is only necessary to describe the differential 
$d_{G,m}^{\,p}: Q(\alpha,p) \ra Q(\alpha',p+1)$ on individual 
summands. This map is zero unless $\alpha' = \alpha$ or $\alpha-1$. 

If $\alpha'=\alpha$, it suffices to describe the map 
\[ \begin{aligned}
{} & \Sym^{m+p+\alpha-n-1}(\Sym^d \,V^*) 
\otimes \wedge^{n+2-p}(\Sym^d \,V^*) \ra \\
& \Sym^{m+p+\alpha-n}(\Sym^d \,V^*) 
\otimes \wedge^{n+1-p}(\Sym^d \,V^*). \end{aligned} \]
This is given by contraction with the canonical trace element in 
\newline $\Sym^d \,V^* \otimes \Sym^d \,V$.
If $\alpha'=\alpha-1$, it suffices to describe the map 
\[ M(\alpha) \otimes \wedge^{n+2-p}(\Sym^d \,V^*) \ra
  M(\alpha-1) \otimes \wedge^{n+1-p}(\Sym^d \,V^*). \]
This is given by contraction with the element in 
$\otimeslim\limits_r \Sym^r(\Sym^{e_r} \,V^*) 
\otimes \Sym^d \,V$, coming from the comultiplication map 
(\ref{comultiplication}).

\subsection{An Alternate Description of $\I_X$.}
We pick up the thread at the beginning of \S \ref{encomplex.subsection}.
There is a bundle map 
\begin{equation} {\widetilde \varphi}: \pi^* (\Omega^1_{\P^d}(1)) \lra 
\ll \label{defn.phitilde} \end{equation}
to be defined as follows. (Let ${\mathfrak T}_\P$ denote the tangent bundle.)
Specifying ${\widetilde \varphi}$ is tantamount to 
specifying a global section in 
\[ H^0(T, \ll \otimes \pi^* {\mathfrak T}_{\P^d}(-1))= 
\otimeslim\limits_{r=1}^d \Sym^r(\Sym^{e_r} \,V^*) 
\otimes \Sym^d \,V, \]  which is tantamount to giving a map 
$\Sym^d \,V^* \lra \otimeslim\limits_{r=1}^d \Sym^r(\Sym^{e_r} \,V^*)$. 
Take this to be the comultiplication map (\ref{comultiplication}).

Over the point $t=(G_1,\dots,G_d;F) \in T$, the transpose map 
$\ll^{-1} \lra \pi^* ({\mathfrak T}_{\P^d}(-1))$ is given by 
\[ {\widetilde \varphi}^{\,\text{trans}}:
\prod\limits_r G_r{}^r \lra \prod\limits_r G_r{}^r +(F) \in 
\frac{\Sym^d \,V}{(F)}. \]
The image is zero iff $F=\prod\limits_r G_r{}^r$ upto constants, i.e., 
iff $t \in \Gamma$. Hence 
the scheme defined by $\,\text{Fitt}_0(\coker\;{\widetilde \varphi})$
is $\Gamma$. 
Now starting with the Koszul complex of ${\widetilde \varphi}$,
the arguments in \S \ref{encomplex.subsection}, \ref{hyperg.subsection} 
go through verbatim. The outcome is a complex of bundles 
${\widetilde \G}^\bullet$ living on $\P^d$, with 
\begin{equation} \begin{aligned}
{\widetilde \G}^p=  
\otimeslim_r \Sym^{r(n+1-p)-e_r-1}& (\Sym^{e_r} V)  \otimes 
\Omega^{n+1-p}_{\P^d} (n+1-p), \\
& \text{\; for $n+1 -d \le p \le n+1 - {\scriptstyle M}$.}
\end{aligned} \label{defn.gtilde} \end{equation}
Then Theorem \ref{hyperg.theorem} and Corollary 
\ref{corollary.gradedpiece} are valid with 
${\widetilde \G}^\bullet$ in place of $\G^\bullet$.

\begin{Remark} \rm 
Consider the map 
\[\gamma_m:\ker\; d_{G,m}^{\,0} \lra  H^0(\P^d,\I_X(m))\]
given by Theorem \ref{hyperg.theorem}.
The inclusion $\ker\; d_{G,m}^{\,0} \subseteq H^0(\P^d, \G^0(m))$ 
is $SL(V)$--equivariant, hence it must split (non-canonically). 
Hence $\gamma_m$ extends to a map 
\[{\widetilde \gamma_m}:H^0(\P^d, \G^0(m)) \lra H^0(\P^d,\I_X(m)).\] 
I believe that it would be of interest to have at least 
one explicit construction of ${\widetilde \gamma_m}$.
\end{Remark}

\section{Some Computations} \label{computation.section}
So far, we have described two approaches to the calculation 
$\I_X$, one in \S \ref{machine.subsection} and the 
other in Theorem \ref{hyperg.theorem}. They can be combined to 
get a concrete invariant--theoretic description of this ideal,
provided  we can get a handle on $\D$. 

We will carry  out the calculation in detail for 
$\lambda=(3\,2)$, and then in addition give the 
result for $(3\,3)$. In each case we will describe the minimal 
resolution of $X_\lambda$ in terms of $SL_2({\mathbb C})$--representations. 
The theorem below will be used for $(3\,2)$, but of course 
it is of more general application. 

\begin{Theorem} \label{bipartite.d.theorem}
Let $\lambda=(\lambda_1,\lambda_2)$ be a partition of $d$, with 
$\lambda_1 \neq \lambda_2$. Then $\D_\lambda$ is supported on the 
rational normal curve $X_{(d)}$, and we have an isomorphism 
$\D_\lambda = {f_{(d)}}_* \O_{\P^1}(-2)$.
\end{Theorem}
\demo 
The claim about supp$(\D)$ will follow directly from Theorem 
\ref{sing.theorem} in the next section. To prove the rest, 
consider the commutative triangle
\[ \diagram
\P^1 \rto^(0.4){\delta} \drto_{f_{(d)}} & \P^1 \times \P^1 \dto^{f_\lambda} \\
& \P^d \enddiagram \] 
Here $\delta$ is the diagonal imbedding, let $\Delta$ be its image. 
\smallskip

\noindent Step 1:
\noindent The ideal of $\Delta$ equals $\O_{\P^1 \times \P^1}(-1,-1)$, 
i.e., we have an extension 
\begin{equation} 
0 \lra \O_{\P^1 \times \P^1}(-1,-1) \lra \O_{\P^1 \times \P^1} 
\lra \O_\Delta \lra 0. 
\label{ext.a} \end{equation}
\newpage 
\noindent Consider the commutative ladder 

\[ \begin{CD}
 0 @>>> \I_{X_{(d)}}/\I_{X_\lambda}  @>>> \O_{X_\lambda} 
@>>> \O_{X_{(d)}}  @>>> 0 \\
 @.   @V{1}VV    @V{2}VV    @V{3}VV @. \\
 0 @>>> {f_\lambda}_* \O_{\P^1 \times \P^1}(-1,-1) @>>> 
{f_\lambda}_* \O_{\P^1 \times \P^1} @>>> {f_\lambda}_* \O_\Delta  
@>>> 0 \end{CD} \]

\noindent whose bottom row is derived from (\ref{ext.a}). 
(The functor ${f_\lambda}_*$ 
is exact since $f_\lambda$ is a finite morphism.)
Now the map $3$ is an isomorphism, so by the five lemma, 
$\coker\, 1 \cong \coker\,2 $. This gives an extension 
\begin{equation}
0 \lra \I_{X_{(d)}}/\I_{X_\lambda} \stackrel{1}{\lra}
{f_\lambda}_* \O_{\P^1 \times \P^1}(-1,-1) \lra \D_\lambda \lra 0.
\label{ext.b} \end{equation}

\noindent Step 2:
\noindent There is an isomorphism 
$\delta^* \O_{\P^1 \times \P^1}(-1,-1) \cong \O_{\P^1}(-2)$, hence by 
taking adjoints, a surjection $\O_{\P^1 \times \P^1}(-1,-1) \ra 
\delta_* \O_{\P^1}(-2)$. The kernel of this surjection equals 
$\O_{\P^1 \times \P^1}(-2,-2)$, so we have an exact sequence 
\[ 0 \lra \O_{\P^1 \times \P^1}(-2,-2) \lra \O_{\P^1 \times \P^1}(-1,-1) \lra 
\delta_* \O_{\P^1}(-2) \lra 0. \]
After applying ${f_\lambda}_*$, we have an extension 
\begin{equation}
 0 \ra {f_\lambda}_* \O_{\P^1 \times \P^1}(-2,-2) \ra 
{f_\lambda}_* \O_{\P^1 \times \P^1}(-1,-1) \ra 
{f_{(d)}}_* \O_{\P^1}(-2) \ra 0. 
\label{ext.c} \end{equation}

\noindent Step 3:
\noindent From (\ref{ext.b}) and (\ref{ext.c}), it suffices to construct an 
isomorphism $\I_{X_{(d)}}/\I_{X_\lambda}
\cong {f_\lambda}_* \O_{\P^1 \times \P^1}(-2,-2)$ 
to prove the theorem. It will be clear from the construction 
below that the isomorphism is compatible with the maps in 
(\ref{ext.b}) and (\ref{ext.c}).

Let $[x_0,x_1; y_0,y_1]$ (respectively 
$[a_0,\dots,a_d]$) be c{\"o}ordinates on $\P^1 \times \P^1$ 
(respectively on $\P^d$). The morphism $f_\lambda$ is defined 
via the identity 
\[ \sum\limits_{j=0}^d a_j \,t^{d-j} = 
(x_0\, t + x_1)^{\lambda_1}(y_0\, t + y_1)^{\lambda_2}.
\] 
We pass to the open set $U_0= \{ a_0 \neq 0 \}$ and write 
\[ \xi= x_1/x_0,\, \eta= y_1/y_0,\;\; \alpha_i=a_i/a_0 \;\;\text{for 
$1 \le i \le d$.} \] The commutative triangle corresponds to a 
diagram of ring homomorphisms 
\[ \diagram
\frac{{\mathbb C}[\xi,\eta]}{(\xi-\eta)} & 
{\mathbb C}[\xi,\eta] \lto_{\text{epi}} \\
& {\mathbb C}[\alpha_1,\dots,\alpha_d] \uto_{q_\lambda} 
\ulto^{q_{(d)}} \enddiagram \]
Let $I_\lambda= \ker \; q_\lambda$, 
$I_{(d)}= \ker \; q_{(d)}$. 

\noindent Claim: The map $q_\lambda: I_{(d)} \lra {\mathbb C}[\xi,\eta]$ has 
image $(\xi- \eta)^2 \,{\mathbb C}[\xi,\eta]$. Hence it defines 
an isomorphism $I_{(d)}/I_\lambda \cong
(\xi- \eta)^2 \,{\mathbb C}[\xi,\eta]$ of 
${\mathbb C}[\alpha_1,\dots,\alpha_d]$--modules.

\noindent Proof:
The ideal $I_{(d)}$ has generators 
\[ J_r = (\alpha_1/d)^r - {\alpha_r}/\tbinom{d}{r}
\quad \text{for\, $2 \le r \le d$.}\]
Set $w=\xi-\eta$, then by straightforward computation,
\[ \begin{aligned}
q_\lambda(\alpha_1^r) & \equiv (d\,\eta)^r + r \lambda_1(d\,\eta)^{r-1}w,\\
q_\lambda(\alpha_r) & \equiv  \tbinom{d}{r}\eta^r + 
\lambda_1\eta^{r-1}\tbinom{d-1}{r-1}w \quad \text{(mod $w^2$)}; 
\end{aligned} \]
hence $q_\lambda(J_r) \equiv 0$\,(mod $w^2$).
We will show that each term $w^2\, \xi^i\eta^j$ lies in 
$q_\lambda(I_{(d)})$. Since 
$\lambda_2 \,\eta= q_\lambda(\alpha_1)-\lambda_1 \,\xi$, it suffices to 
show this for $j=0$. Due to the equality 
\[ (\lambda_1^2+\lambda_1\lambda_2)\, \xi^2 = 
2\lambda_1\,q_\lambda(\alpha_1)\,\xi + (\lambda_2-1) \,q_\lambda(\alpha_1^2)
-2\lambda_2\,q_\lambda(\alpha_2), \] we may further assume $i=0,1$.
Now a tedious calculation shows that a constant multiple 
of $q_\lambda(J_2)$ equals $w^2$, and 
a certain linear combination of 
$q_\lambda(J_3)$ and $q_\lambda(\alpha_1\,J_2)$ equals $w^2\,\xi$.
This proves the claim. 

Evidently, the same calculation goes through over the open 
set $U_d= \{a_d \neq 0 \}$, where we have an isomorphism 
\[ I_{(d)}/I_\lambda \cong
(1/\xi- 1/\eta)^2 \,{\mathbb C}[1/\xi,1/\eta].\]
The union $U_0 \cup U_d$ misses only two 
points of $X_\lambda$; say $P$ and $Q$, which are the images of points 
$[1,0] \times [0,1],[0,1] \times [1,0] \in \P^1 \times \P^1$.
We have shown that over the open set $X_\lambda - \{P,Q\}$, 
we have an isomorphism 
$\I_{X_{(d)}}/\I_{X_\lambda}={f_\lambda}_* \ll$
for some line bundle $\ll$. 

If $g$ is an automorphism of $X_\lambda$ coming from 
$SL(V)$, then $g^*(\I_{X_{(d)}}/\I_{X_\lambda})$ $\cong 
\I_{X_{(d)}}/\I_{X_\lambda}$. We can move $P$ or $Q$ 
to a point of $U_0 \cup U_d$ by such a $g$, so the isomorphism 
is valid throughout $X_\lambda$. 
Now take the ratio of the local generators 
\[ \frac{(x_1/x_0 - y_1/y_0)^2}{(x_0/x_1 - y_0/y_1)^2} 
= (x_0/x_1)^{-2}(y_0/y_1)^{-2}, \]
which shows that $\ll=\O_{\P^1 \times \P^1}(-2,-2)$.
The map 1 in (\ref{ext.b}) is the restriction of map 2, and 
locally, the latter is just $q_\lambda$ at the level of rings. 
Hence the compatibility claimed earlier. This completes the 
proof of the theorem. \qed 
\begin{Corollary} 
For $\lambda$ as above, 
$H^0(\P^d, \D_\lambda(m)) = \Sym^{dm-2} \,V^*$ for 
any $m \ge 0$. The Hilbert polynomial of $X_\lambda$ is 
$\lambda_1 \lambda_2\, m^2 + 2$. 
\end{Corollary}
\demo 
We have $\D_\lambda \otimes \O_{\P^d}(m) =
{f_{(d)}}_* \O_{\P^1}(dm-2)$ from the projection formula, 
hence the first claim. Now the second claim is a 
simple corollary of formula (\ref{defn.dlambda}). \qed

\medskip 
The following device will be useful for calculating with 
complexes of representations\,:
\medskip 

\noindent Let ${\mathbb A}$ be the free abelian group on symbols 
$\{s_n \}_{n \ge 0}$. We let $\chi(\Sym^n \,V)=s_n$ for any $n$, and then 
define $\chi(R) \in {\mathbb A}$ for any finite dimensional 
$SL(V)$--representation $R$ by decomposing it into irreducibles. 
Thus $\chi(R)$ is merely the formal character of $R$.

If $G^\bullet$ is a bounded complex of $SL(V)$--representations and 
$SL(V)$--equivariant cochain maps, then define 
$\chi(G^\bullet)= \sum\limits_p (-1)^p \chi(G^p)$. This can be seen 
as the `Euler characteristic' of $G^\bullet$, taking values in 
${\mathbb A}$. The following equality is immediate 
from Schur's lemma\,:
\begin{equation} \label{wordequality}
\chi(G^\bullet) = \sum\limits_p \,(-1)^p \chi(H^p(G^\bullet)). 
\end{equation}

\subsection{The Surface $X_{32}$.} \label{surface32.calc}
As in \S \ref{machine.subsection}, 
we have a morphism 
\[  f_{32}: \P\, V \times \P\, V {\lra} 
   \P\, \Sym^5 \,V, \quad (G_3,G_2) \lra G_3{}^3.G_2{}^2 \,(=F); \] 
and in affine c{\"o}ordinates, a ring map 
\[  q_{32}: {\mathbb C}[a_1,\dots,a_5] \lra
{\mathbb C}[u,v], \]
\[ a_1 \ra 3u+2v,\, a_2 \ra 3u^2+6uv+v^2,\dots,\, a_5 \ra u^3v^2. 
\] 
Homogenize the kernel of $q_{32}$ with respect to $a_0$, and 
take its minimal resolution:
\[ \begin{aligned} 
0 & \lla \I_X \lla \O_\P(-4)^{28} \lla
\O_\P(-5)^{68}  \lla \O_\P(-6)^{64} \lla 
\O_\P(-7)^{28} \\ & \lla \O_\P(-8)^5  \lla 0. 
\end{aligned} \]
All of this was carried out in Macaulay--2.
We write this as 
\[ 0 \la \I_X \la \E^0 \la \E^{-1} \dots \la 
\E^{-4} \la 0. \]
Our task now is to identify the $SL(V)$--representations 
entering into the syzygy modules $\E^p$. There is a 
spectral sequence 
\begin{equation}\label{spec.calc}
 \begin{aligned} E_1^{p,q} & = H^q(\P^5,\E^p(m)),\quad 
    d_r^{\,p,q}: E_r^{p,q} \lra E_r^{p+r,q-r+1} \\
E_\infty^{p,q} & \RA H^{p+q}(\P^5,\I_X(m)).
\end{aligned} \end{equation}
Here $\E^p(m)$ denotes $\E^p\otimes \O_{\P^5}(m)$ etc.
Note the equality
\begin{equation} \begin{aligned}
{} & \chi(H^0(\G^\bullet(m))) \;
(=\sum\limits_p \,(-1)^p \chi(H^0(\P,\G^p(m)))\;\;\text{by definition}) \\
=\, & \chi(\Hy^0(\G^\bullet(m))) - 
\chi(\Hy^1(\G^\bullet(m))),
\end{aligned} \label{wordequalityforG} \end{equation}
which is merely (\ref{wordequality}) rewritten using 
Theorem \ref{hyperg.theorem}.

Now let $m=4$. From (\ref{spec.calc}) we have 
$H^{\ge 1}(\I_X(4)) =0$, so $H^1(\O_X(4))=H^2(\I_X(4))=0$. 
(We are using nothing more than the standard computation of 
cohomology of line bundles on projective spaces, 
see \cite[\S 3.5]{Ha}.) Now by 
Theorem \ref{hyperg.theorem},
\[ \chi(\Hy^1(\G^\bullet(4)))= \chi(H^0(\P,\D(4))) =
\chi(\Sym^{18} \,V^*) =s_{18}. \]
The alternating sum $\chi(H^0(\G^\bullet(4)))$ can be calculated 
directly from (\ref{qalphap.formula}) and repeated use of 
formulae (\ref{slv.identities}) and (\ref{plethysm}). 
(This was programmed in Maple.)
Substituting this into (\ref{wordequalityforG}), we get 
\[ \chi(H^0(\I_X(4)))= \chi(\Hy^0(\G^\bullet(4))) 
= s_{12}+s_8+s_4+s_0. \] 

Now choose $m=5$. From (\ref{spec.calc}),
we deduce $H^{\ge 1}(\I_X(5))=0$ and an extension 
\[  0 \lra H^0(\E^{-1}(5)) \lra H^0(\E^0(5)) 
\lra H^0(\I_X(5)) \lra 0. \] 
By the previous step, 
\[ H^0(\E^0(5)) = \Sym^5 \,V^* \otimes 
(\Sym^{12} \,V^* \oplus \Sym^8 \,V^* \oplus \Sym^4 \,V^* 
\oplus {\mathbb C}). \]

Now $\chi(H^0(\G^\bullet(5)))$ and $\chi(\Hy^1)$ are calculated 
exactly as before, and one can identify the unknown term 
\[ \chi(H^0(\E^{-1}(5)))= 
s_{13}+s_{11}+s_9+2s_7+2s_5+s_3. \]
It is clear how to proceed to determine the entire minimal 
resolution. The complete result is as follows.
(We have written e.g. $\{5,3^2,1^3\}$ for  
the module $\Sym^5 \,V^* \oplus (\Sym^3 \,V^*)^{\oplus 2} 
\oplus (V^*)^{\oplus 3}$.)
\[ \begin{aligned} 0 & \lla \I_X \lla  \{12,8,4,0 \} \otimes \O_\P(-4) \lla
\{ 13,11,9,7^2,5^2,3 \} \otimes \O_\P(-5) \\ & \lla 
  \{ 12,10,8^2,6,4^2,2,0^2 \} \otimes \O_\P(-6) \lla 
\{ 9, 7,5,3 \} \otimes \O_\P(-7) \\ & \lla 
\{ 4 \} \otimes \O_\P(-8) \lla 0. 
\end{aligned} \]
It follows then, that $I_{X_{32}}$ is generated by a 
subrepresentation 
\[ (I_X)_4=
\Sym^{12} \oplus \Sym^8 \oplus \Sym^4 \oplus \Sym^0 
\subseteq \Sym^4(\Sym^5 \,V^*). \]
We can identify the summands in the language of classical 
invariant theory. Define covariants 
\[ H=(F,F)^2, \;\; i=(F,F)^4, \;\; A=(i,i)^2 \]
for the generic binary quintic $F$.
(In classical literature the notation varies from source to source. 
I have consistently followed \cite[Ch. VII]{GrYo}, except that their $f$ is 
our $F$. All the transvectants are calculated using the symbolic method, see 
Scholium \ref{transvectant.remark}.) 

The inclusion $\Sym^{12} \subseteq \Sym^4(\Sym^5 \,V^*)$ 
corresponds to a covariant of degree $4$ and order $12$, 
we tentatively call it $J$. (Said differently, 
let $J= \sum_i \varphi_i\,x^iy^{12-i}$. Then the vector space 
spanned by the forms $\varphi_i$ is a subrepresentation of $(I_X)_4$,
isomorphic to $\Sym^{12}\, V^*$.) 
Now the space of covariants of  
type $(4,12)$ is spanned by $H^2$ and $i.F^2$ (loc.cit.), hence 
$J= \alpha\, H^2 + \beta\, i.F^2$ for some constants 
$\alpha,\beta$. Choose a form $F$ lying in $X_{32}$
(but otherwise sufficiently general) and evaluate 
the right hand side. By hypothesis $J$ must vanish, from this 
we can deduce that $\alpha:\beta= 25:-6$. Thus the inclusion 
$\Sym^{12} \subseteq \Sym^4(\Sym^5 \,V^*)$ is entirely specified 
by the covariant $25 \,H^2 - 6\,i.F^2$. The other summands in 
$(I_X)_4$ are calculated in the same way, and we get the 
following theorem:
\begin{Theorem} 
A homogeneous binary quintic $F$ can be factored as 
$L^3 M^2$ for some linear forms $L,M$; if and only if, 
the covariants 
\[ 25\,H^2-6\,i.F^2, \; 5\,i.H+6\,F.(i,F)^2, \;
2\,i^2+15\,(i,H)^2 \;\; \text{and $A$} \]
vanish on $F$.
\label{criterion.32} \end{Theorem}
These are the `algebraic conditions' on $F$, alluded to in the 
very beginning of the paper. 

\subsection{The surface $X_{33}$.}
We will merely state the result for the surface $X_{33}$, 
suppressing all computation. (This case is in fact easier than 
the previous one, since $\D_{33}=0$.) The equivariant 
minimal resolution turns out to be 
\[ \begin{aligned} 
{} & 0 \la \I_X \la 
\{ 12,8,6 \} \otimes \O_\P(-3) \la 
\{14,12,10^2,8,6^2,4,2^2 \} \otimes \O_\P(-4) \\ & \la 
\{14,12,10^2,8^2,6^3,4^2,2\} \otimes \O_\P(-5) \la 
\{12,10,8^2,6,4^2,2,0^2 \} \otimes \O_\P(-6) \\ & \la
\{8,6,4\} \otimes \O_\P(-7) \la 
\{2\} \otimes \O_\P(-8) \la 0.
\end{aligned} \]
Thus the ideal is generated in degree three by three 
covariants, and these can be identified as before.
Again let $H=(F,F)^2$ and $i=(F,F)^4$ for the generic 
binary sextic $F$.
\begin{Theorem} 
A homogeneous binary sextic $F$ can be written as 
the cube of a quadratic form, if and only if, 
the covariants 
\[ (F,H),\, (F,i)\, \;\text{and $8\,F.(F,F)^6 -75\,(F,i)^2$} \]
vanish on $F$.
\label{criterion.33} \end{Theorem}

It is almost superfluous to mention that Theorems \ref{criterion.32}
and \ref{criterion.33} are not being seriously suggested as
\emph{practical criteria} for detecting multiple factors.

\begin{Remark} \rm 
Here we give the result for the case $\lambda=(3\,2\,2)$,
which was left unfinished in \S \ref{machine.subsection}.
The sheaf $\D_{322}$ sits in an extension 
\[ 0 \ra {f_{52}}_* \O_{\P^1 \times \P^1}(-3,-1) 
\ra \D_{322} \ra {f_7}_* \O_{\P^1}(-2) \ra 0. \]
This is proved by an argument similar to (but more 
elaborate than) the one used in Theorem 
\ref{bipartite.d.theorem}. From this, we deduce an 
exact sequence 
\[ 0 \ra \Sym^{5m-3} \,V^* \otimes \Sym^{2m-1} \,V^* \ra 
H^0(\P^7, \D(m)) \ra 
\Sym^{7m-2} \,V^* \ra 0, \] and now the calculation proceeds 
as before. 
The ideal is generated by the graded piece 
\[ \begin{aligned} 
(I_X)_6 & = 
\{ 30, 26,24,22^2,20,18^3,16,14^3,12^2,10^3,8,6^3,2^2 \} \\
& \subseteq \Sym^6(\Sym^7 \,V^*). \end{aligned} \]
\end{Remark} \rm 
\begin{Remark} \rm 
In the case dim $X_\lambda =1$ (i.e., of the rational 
normal curve) the minimal resolution is very pretty, so I
cannot resist recording it here. Its ideal is generated by 
quadrics, more specifically 
\[ (I_X)_2= \opluslim\limits_{r=1}^{[d/2]} \Sym^{2d-4r} \,V^* 
\subseteq \Sym^2(\Sym^d \,V^*).\] 
(The summands correspond to the transvectants $(F,F)^{2r}$ for 
$1 \le r \le [d/2]$.)
The natural multiplication map 
\[\Sym^{d-1} \,V^* \otimes V^* \lra \Sym^d \,V^* =
H^0(\P^d,\O_\P(1)),\]
gives rise to a bundle map 
\[\varphi: \O_\P(-1) \otimes \Sym^{d-1} \,V^* \lra V \otimes \O_\P. \]
Then the inequality rank $\varphi \le 1$ holds over a point 
$[w] \in \P\, \Sym^d\,V$, if and only 
if, $w=v^d$ for some $v \in V$. Hence 
the Eagon--Northcott complex of $\varphi$ resolves $\O_X$.
\end{Remark}

I enclose a list of representations entering into the 
minimal generators of $I_X$ and the corresponding covariants 
when $\lambda$ has two parts and 
$4 \le d \le 6$. For each $d$, define 
$H=(F,F)^2,\,i=(F,F)^4$. In each case $I_m$ denotes the degree 
$m$ part of $I_X$, the rest needs no explanation. 
\[ \begin{aligned}
\text{Case $d=4$.}&  \\
(31): \quad I_2=& \;\{0\} = \{ i\}, \,I_3=\{0\} = \{ (F,H)^4 \}. \\
(22): \quad I_3=& \;\{6\}=\{(F,H)\}. \\ \\
\text{Case $d=5$.}& \\
(41): \quad I_2=& \;\{2\} = \{ i\}. \\
(32): \quad I_4=& \;\{12,8,4,0\} \\
=& \;\{ 25\,H^2-6\,i.F^2, 5\,i.H+6\,F.(i,F)^2, 
2\,i^2+15\,(i,H)^2, A \}. \end{aligned} \]
\[ \begin{aligned} 
\text{Case $d=6$.}&  \\
(51): \quad I_2=& \;\{4,0\} = \{ i, (F,F)^6\}. \\
(42): \quad I_2=& \;\{0\}=\{(F,F)^6\},\;
            I_3= \;\{2\}=\{(F,i)^4\}, \\
            I_4=& \;\{16,12,0\}=
                   \{27\,H^2-8\,i.F^2,3\,i.H+4\,F.(F,i)^2, (i,i)^4 \}. \\
(33): \quad I_3=& \;\{12,8,6\} 
               = \{ (F,H), (F,i), 8\,F.(F,F)^6 -75\,(F,i)^2 \}. 
\end{aligned} \]

\subsection{The Castelnuovo regularity of $\I_X$.}
\label{cmreg.subsection}

Let $\F$ be a coherent $\O_\P$--module. It is said to be $m$-regular
for an integer $m$, if $H^i(\P,\F(m-i))=0$ for $i \ge 1$. If 
this is so, then $\F(m)$ is generated by global sections and $\F$ is 
$m'$-regular for any $m' \ge m$. 
This concept is due to David Mumford, see \cite{Mum} 
for a masterly exposition. 

Here we state a result bounding the regularity of $\I_{X_\lambda}$, when 
$\lambda=(\lambda_1,\lambda_2)$. Define 
\[ \mu = \begin{cases} 
d{}^2/4-d+3 & \quad \text{for $\lambda_1 = \lambda_2$,} \\
\lambda_1 \lambda_2 (2\lambda_1 \lambda_2-d+2)^2 + 
2\lambda_1 \lambda_2 -d+5 & \quad \text{otherwise.} 
\end{cases} \]
\smallskip

\begin{Theorem} With $\lambda$ as above, $\I_X$ is $\mu$-regular. 
A fortiori, $I_X$ is generated by forms of degree 
$\le \mu$. \end{Theorem}
\demo The first case follows directly from the central result of 
\cite{Lazarsfeld1}. For the second case, note the following lemma.

\begin{Lemma}[\cite{KL.SGA6},\cite{Mum}]\label{reg.bound.lemma}
Let $\F$ be a coherent $\O_\P$-module and let $H$ be a general 
hyperplane in $\P$. (It is sufficient that the generic point of 
$H$ not be associated to $\F$.) Assume that $\F|_H$ is $m$--regular.
Then $H^{\ge 2}(\P,\F(m-1))=0$ and $\F$ is 
$(h^1(\F(m-1))+m)$--regular. \qed \end{Lemma} 

Now $\I_X|_H$ is the ideal of a nondegenerate curve of degree 
$2 \lambda_1 \lambda_2$ in $\P^{d-1}$; hence by 
\cite[Theorem 1.1]{GPL} it is
$m_0$--regular, where $m_0=(2 \lambda_1 \lambda_2-d+3)$. 
From the surjection $H^0(\O_X(m_0-1)) \ra H^1(\I_X(m_0-1)) \ra 0$, 
we have $h^1(\I_X(m_0-1)) \le h^0( \O_X(m_0-1))$. 
Since $H^{\ge 2}(\I_X(m_0-1))=0$, we deduce 
\[
h^0(\O_X(m_0-1)) = \chi(\O_X(m_0-1)) = 
\lambda_1\lambda_2(m_0-1)^2 +2. \]
The theorem is proved.  \qed 
\medskip 

It is not likely that the bound is optimal, see the last 
section for a conjecture. For any $\lambda$, if the Hilbert polynomial of 
$\D_\lambda$ is known, then such a bound can be formulated. See 
\cite{KL.SGA6} for a general regularity bound along such 
lines\footnote{However, the definition of 
`$(b)$-polyn{\^o}me' and the statement of `th{\'e}or{\`e}me 
principal' are not very clear to me as stated there.}.

\section{Singularities of the CR locus}
\label{singularities.section}
In case of a general partition $\lambda$, an explicit determination 
of $\D_\lambda$ on the lines of Theorem \ref{bipartite.d.theorem} 
does not appear to be easy. (The isomorphism 
$\H^1(\G^\bullet)=\D$ is of little help in this.) In this section, we 
solve the easier problem of describing the support of 
$\D_\lambda$,
which coincides with the singular locus of $X_\lambda$.

Note that a point $F \in X_\lambda$ is nonsingular iff $\O_{X,F}$ is 
normal. 

\begin{Proposition} \label{sing.prop0}
A point $F \in X$ is nonsingular if and only if the following 
conditions hold: 
\begin{enumerate}
\item{The preimage $f_\lambda^{-1}(F)$ is a singleton set, 
say $\{ G \}$; and }
\item{the map $df_\lambda: T_{Y,G} \lra T_{\P^d,F}$ on 
tangent spaces is injective.}
\end{enumerate}\end{Proposition} 
\demo This is immediate from the local criterion of isomorphism 
in \cite[Theorem 14.9]{JoeH}. \qed 

With the partition $\lambda=(1^{\,e_1} 2^{\,e_2}\dots d^{\,e_d})$, 
we now associate a set $\ss_\lambda$ of partitions of $d$. The point 
of the definition lies in Theorem \ref{sing.theorem}.

\begin{Definition} \rm 
A partition $\mu=(1^{f_1} 2^{\,f_2} \dots d^{\,f_d})$ belongs to 
$\ss_\lambda$ in one of the following (mutually exclusive) cases.
\begin{enumerate}
\item[a.]\label{defmerger.a}
{ There exist distinct integers $r_1,r_2$ such that 
$\mu$ is derived from $\lambda$ by setting 
\[ \begin{array}{llll}
f_{r_1} & =  e_{r_1}-1, \quad & f_{r_2} & =  e_{r_2}-1, \\ 
f_{r_1+r_2} & =  e_{r_1+r_2}+1 \;\;\text{and}  &  
f_r & =  e_r\; \text{ elsewhere.} \end{array} \] }
\item[b.]\label{refmerger.b}
{ There exist integers $r_1,r_2,t$\, such that 
$r_1 \neq r_2, \,e_{r_1} > 0$ and $\mu$ is derived from 
$\lambda$ by setting 
\[ \begin{array}{llll}
f_{r_1} & =  e_{r_1}+1, \quad & f_{r_2} & =  e_{r_2}-t 
\;\;\text{and} \\
f_r & =  e_r\; \text{ elsewhere.} \end{array} \]
(This forces $r_1=t\,r_2$.) }
\item[c.]\label{refmerger.c}
{ There exist integers $r_1,r_2,r_3,t_1,t_2$ such that 
$r_1,r_2,r_3$ are pairwise distinct, 
$r_3=t_1\,r_1=t_2\,r_2$, and $\mu$ is derived from $\lambda$ by setting 
\[ \begin{array}{llll}
f_{r_1} & =  e_{r_1}-t_1, \quad & f_{r_2} & =  e_{r_2}-t_2, \\
f_{r_3} & =  e_{r_3}+2 \;\;\text{and} & 
f_r & =  e_r\; \text{ elsewhere.} \end{array} \] }
\end{enumerate}
\end{Definition}
In each case $\lambda$ is a refinement of $\mu$, so 
$X_\mu \subseteq X_\lambda$. We will write 
$\ss_\lambda =\ss_\lambda^{(a)} \cup \ss_\lambda^{(b)} 
\cup \ss_\lambda^{(c)}$. Then $\ss_\lambda$ is empty if and only 
if all parts in $\lambda$ are equal.

\begin{Example} \rm 
Let $\lambda=(1\, 2^3\, 3^2\, 4)$. Then $\ss_\lambda^{(a)}$ has in all 
six elements, $(1\, 2^3\, 3\, 7)$,\,$(1\, 2^2\, 3\, 4\, 5)$ etc.
We can set $f_2=1,f_4=2$ and $f_r=e_r$ elsewhere, hence 
$(1\, 2^2\, 3^2\, 4^2) \in \ss_\lambda^{(b)}$. We can set 
$f_6=2,f_3=f_2=0$ and $f_r=e_r$ elsewhere, hence 
$(1\,4\,6^2) \in \ss_\lambda^{(c)}$.

\end{Example} 

\begin{Theorem} \label{sing.theorem}
The singular locus of $X_\lambda$ equals \, $\bigcup\limits_\mu X_\mu$, 
\;the union quantified over all $\mu \in \ss_\lambda$. 
\end{Theorem}
The union may well be redundant, i.e., there may exist $\mu_1,\mu_2$
such that $X_{\mu_1} \subseteq X_{\mu_2}$.
The theorem is a consequence of the following two propositions:
\begin{Proposition} \label{sing.prop1}
Let $F$ be a singular point of $X_\lambda$. 
\begin{enumerate}
\item{If $f_\lambda^{-1}(F)$ is 
singleton, then $F \in X_\mu$ for some $\mu \in \ss_\lambda^{(a)}$.}
\item{ If $f_\lambda^{-1}(F)$ has more than one element and 
$F \notin X_\mu$ for any  $\mu \in \ss_\lambda^{(a)}$, then 
$F \in X_\mu$    for some $\mu \in \ss_\lambda^{(b)} \cup \ss_\lambda^{(c)}$.}
\end{enumerate} \end{Proposition}

\begin{Proposition} \label{sing.prop2}
Let $F \in \bigcup\limits_{\mu \in \ss_\lambda} X_\mu$ be a point. 
Then either $f_\lambda^{-1}(F)$ is not singleton or 
$T_{Y,G} \lra T_{\P^d,F}$ is not injective for the unique 
$G$ in $f_\lambda^{-1}(F)$. Hence in either case 
$F$ is a singular point of $X$.
\end{Proposition}

\medskip 
The following lemma will be needed in the proofs.
\begin{Lemma} 
Let $\P\, \Sym^m \,V \times \P\, \Sym^n \,V 
\stackrel{h}{\lra} \P\, \Sym^{m+n} \,V$ denote the 
multiplication map $(A,B) \lra A.B(=C)$. The map 
\[dh: T_{\P^m \times \P^n,(A,B)} \lra T_{\P^{m+n},C} \]
on tangent spaces is injective, if and only if $A,B$ have 
no common factor.
\end{Lemma}
\demo The map is given as 
\[ ([\alpha_0,\dots,\alpha_m],[\beta_0,\dots,\beta_n]) \lra 
[\gamma_0,\dots,\gamma_{m+n}], 
\quad \gamma_k=\sum\limits_{i+j=k} \alpha_i \beta_j. \]
By a linear change in $x,y$ we may assume that $\alpha_0=\beta_0=1$ 
at $(A,B)$. 
Then ${\partial \gamma_k}/{\partial \alpha_i}=\beta_{k-i}$,
${\partial \gamma_k}/{\partial \beta_j}=\alpha_{k-j}$, and the Jacobian matrix 
representing $dh$ at $(A,B)$ equals the Sylvester resultant of $A$ and $B$.
Hence it has rank $m+n$ exactly when $A,B$ have no common factor.
\qed 

\begin{Corollary} \label{sing.cor} 
If $G=(G_1,\dots,G_d) \in Y_\lambda$ is a point with $f_\lambda(G)=F$; 
then $df:T_{Y,G} \lra T_{\P^d,F}$ is injective if and only if 
no two of the $G_r$ have a common factor. 
\end{Corollary}
\demo The map $\Pi v_r$ (see (\ref{defn.nulambda})) is injective 
on tangent spaces everywhere. Now write $\mu$ as a succession of maps 
with products of two factors at a time, and apply the lemma repeatedly.
\qed

\medskip 
In the sequel, we will use the same notation for 
two forms which differ by a nonzero multiplicative constant. Since there will
be no occassion to add two forms, this should cause no confusion.
If $G$ is a form and $L$ a linear form, then $\ord_L(G)=m$ means 
that $L^m \vert \,G$, but $L^{m+1} \slash \hspace{-1.5mm}  \vert \,G$. 

\medskip 
{\sc Proof of proposition \ref{sing.prop1}.}

1. Say $f^{-1}(F)=\{G\}$, then by Prop. \ref{sing.prop0}
and the previous corollary, two of the $G_r$ have a common factor. 
Say a linear form $L$ divides $G_{r_1}, G_{r_2}$.  Define 
\[ \begin{array}{llll}
Q_{r_1} & =  G_{r_1}/L, \quad & Q_{r_2} & =  G_{r_2}/L, \\ 
Q_{r_1+r_2} & =  L.G_{r_1+r_2} \;\;\text{and}  &  
Q_r & =  G_r\; \text{ elsewhere.} \end{array} \] 
Then $Q$ defines a point of $f_\mu^{-1}(F)$ with $\mu \in \ss_\lambda^{(a)}$.
\qed 
\medskip 

2. Let $G=(G_r),\,H=(H_r)$ be two distinct points in 
$f_\lambda^{-1}(F)$. No two of the $G_r$ have a common 
factor (otherwise $F$ would lie in $X_\mu$ for some 
$\mu \in \ss_\lambda^{(a)}$) and similarly for 
the $H_r$. Let $r_1$ be the largest index such 
that $G_{r_1} \neq H_{r_1}$. (This of course means $G_{r_1}$ is not a 
scalar multiple of $H_{r_1}$.) Then some linear factor 
$L$ of $G_{r_1}$ is not a factor of $H_{r_1}$. (Assume the contrary. 
Then some linear factor, say $K$, 
occurs with a higher power in $G_{r_1}$ than in $H_{r_1}$. But then 
$K$ must divide some $H_{r'}$ for $r' \neq r_1$. This is impossible since 
$H_{r_1},H_{r'}$ have no common factors.) By the same argument, there exists 
a linear form $M$ dividing $H_{r_1}$ but not $G_{r_1}$. 
Say $\ord_L(G_{r_1})=\alpha, \,
\ord_M(H_{r_1})=\beta$ and without loss of generality $\alpha \le \beta$.

Case $\alpha=1$. Now $L$ divides $H_{r'}$ for exactly one 
value of $r'$, let $r_2$ be this value. If $\ord_L(H_{r_2})=t$, then 
$t\,r_2 = r_1$. Define 
\[ \begin{array}{llll}
Q_{r_1} & =  L.H_{r_1}, \quad & Q_{r_2} & =  H_{r_2}/L^t \;\;\text{and} \\ 
Q_r & =  H_r\; \text{ elsewhere.} \end{array} \] 
Then $Q$ is a point of $f_\mu^{-1}(F)$ with $\mu \in \ss_\lambda^{(b)}$.

Case $\alpha >1$. Set $t_1=\alpha$ and $r_3=t_1r_1$.
Assume $\ord_L(H_{r_2})=t_2 >0$ for some $r_2$ as above, then 
necessarily $r_2t_2=r_3$. Let 
\[ \begin{array}{llll}
Q_{r_1} & =  H_{r_1}/M^{t_1}, \quad & Q_{r_2} & =  H_{r_2}/L^{t_2}, \\
Q_{r_3} & =  L.M.H_{r_3} \;\; \text{and} & 
Q_r & =  H_r \quad \text{ elsewhere.} \end{array} \] 
Then $Q$ is a point of $f_\mu^{-1}(F)$ with 
$\mu \in \ss_\lambda^{(c)}$. \qed

\medskip 
{\sc Proof of proposition \ref{sing.prop2}.}
Assume that $F \in X_\mu$ for some $\mu \in \ss_\lambda$ and that 
$f_\lambda^{-1}(F)$ is a singleton set $\{ G \}$.
We will deduce that two of the $G_r$ must have 
a common factor; this will prove the proposition. 
In the sequel, $H$ denotes a point of $f_\mu^{-1}(F)$.

Case \;$\mu \in \ss_\lambda^{(a)}$. 
Choose a linear form $L$ dividing $H_{r_1+r_2}$, and define 
\[ \begin{array}{llll}
Q_{r_1} & =  L.H_{r_1}, \quad & Q_{r_2} & =  L.H_{r_2}, \\
Q_{r_1+r_2} & =  H_{r_1+r_2}/L \;\;\text{and} & 
Q_r & =  H_r \quad \text{ elsewhere.} \end{array} \] 
Then $Q \in f_\lambda^{-1}(F)$, so $Q=G$. But then 
$G_{r_1},G_{r_2}$ have the common factor $L$. 

Case \;$\mu \in \ss_\lambda^{(b)}$. 
Note that $\deg H_{r_1} \ge 2$, and let $L,L'$ be any two linear factors 
of $H_{r_1}$. We claim that $L = L'$. Suppose not, and 
define a point $Q$ by 
\[ \begin{array}{llll}
Q_{r_1} & =  H_{r_1}/L, \quad & Q_{r_2} & =  H_{r_2}.L^t \quad\text{and} \\
Q_r & =  H_r \quad \text{ elsewhere;} \end{array} \] 
and $Q'$ by the same formulae with $L'$ instead of $L$. 
Then $Q,Q' \in f_\lambda^{-1}(F)$ and $Q_{r_1} \neq Q_{r_1}'$,
which is impossible. Hence necessarily $L = L'$, which implies that 
$L^2 | H_{r_1}$. But now $Q=G$, and $Q_{r_1},Q_{r_2}$ have the common 
factor $L$.

Case \; $\mu \in \ss_\lambda^{(c)}$. Let $L_1,L_2$ be two linear factors 
of $H_{r_3}$ and define 
\[ \begin{array}{llll} 
Q_{r_1} & =  H_{r_1}.L_1^{t_1}, \quad & Q_{r_2} & =  H_{r_2}.L_2^{t_2}, \\
Q_{r_3} & =  H_{r_3}/{L_1 L_2} \;\;\text{and} & 
Q_r & =  H_r \quad \text{elsewhere.} \end{array} \] 
Define $Q'$ by the same formulae with $L_1,L_2$ interchanged. Then 
$Q,Q' \in f_\lambda^{-1}(F)$, hence they must both equal $G$. Hence 
$L_1 = L_2$, and $G_{r_1},G_{r_2}$ have the common factor $L_1$.
\qed 

\medskip 
The proof of Theorem \ref{sing.theorem} is complete. 
On intuitive grounds it is plausible that the CR locus should develop a 
singularity wherever some parts in $\lambda$ merge together. 
The element of surprise (for the author) lies in that only certain 
mergers produce singularities.

\section{Comments and Conjectures}
\label{comments.section}
In his memoir \cite{Cayley1}, Cayley gives necessary and sufficient 
conditions for a form to lie in $X_\lambda$, 
for all partitions $\lambda$ of $d \le 5$. This amounts to 
describing $I_X$ upto radicals. (Of course, his paper predates the 
recognition of the distinction between set--theoretic and ideal--theoretic
generation.) E.g., for $\lambda=(3\,2)$, he proves that the 
summand $\Sym^{12} \subseteq \Sym^4(\Sym^5 \,V^*)$ alone defines 
$I_X$ set--theoretically\footnote{Since Cayley's notation is different from 
ours, some effort is required to see that the first covariant in 
Theorem \ref{criterion.32} agrees with his solution.}.
However, \emph{pace} Cayley's own assertion, I do not see how to make 
his method work for a general $\lambda$.

J. Weyman has proved \cite[Theorem 3]{Weyman1} that 
for $\lambda=(r\,1^{d-r})$ with $r \ge [d/2]+1$, the ideal 
$I_{X_\lambda}$ is generated in degrees $\le 4$. (Also 
see \cite{Weyman3} and \cite{Weyman2} for 
results about the generators and Hilbert function of $I_X$.)
I hazard the following conjecture (which I have numerically 
verified for $d \le 8$).
\begin{Conjecture} \rm 
Let $\lambda=(\lambda_1,\lambda_2)$ be a partition of $d$. If 
$\lambda_1 \neq \lambda_2$ (respectively $\lambda_1=\lambda_2$), 
then the sheaf $\I_X(4)$ (resp. $\I_X(3)$) is Castelnuovo regular. 
In particular, $I_X$ is generated in degrees $\,\le 4$ (resp. 
degrees $\,\le 3$).
\end{Conjecture}

It is proved in \cite{Schreyer1}, that the variety $X_{(d-1,\,1)}$
is arithmetically Gorenstein. It is arithmetically Cohen--Macaulay 
in two more simple cases ($\lambda=(d),(2\,1^{d-2})$), but 
does not seem to be so for any other case. (I have numerically 
checked this for $d \le 6$.) One would like to have a proof or 
a counterexample. 

It would be valuable to have a structure theorem for 
$\D_\lambda$ similar to Theorem \ref{bipartite.d.theorem}. 
There are a few more classes of partitions for which this can 
be achieved, and I hope to return to this 
problem in future. 
\bibliographystyle{plain}
\bibliography{../BIB/hema1,../BIB/hema2,../BIB/hema3}
\smallskip

\centerline{------------------------------}

\medskip
\parbox{15cm}{ \small
Department of Mathematics and Statistics, \\ Queen's University,
Kingston, Ontario K7L 3N6. Canada. \\
email: jaydeep@mast.queensu.ca }
\end{document}